\documentclass[12pt,reqno]{amsart}
\usepackage{amsmath}
\usepackage[english,activeacute]{babel}
\usepackage{inputenc}
\usepackage{amssymb}
\usepackage{amsthm}
\usepackage{mathtools}
\usepackage{graphics,graphicx,tikz}
\usepackage{array}
\usepackage{a4wide}
\allowdisplaybreaks
\usepackage{color, url}
\usepackage{float}
\usepackage{multicol}
\usepackage[shortlabels]{enumitem}
\usepackage[colorlinks=true, citecolor=red, linkcolor=blue, urlcolor=blue]{hyperref}
\theoremstyle{plain}
\newtheorem{theorem}{Theorem}[section]

\theoremstyle{definition}

\newtheorem{conjecture}[theorem]{Conjecture}

\begin{document}
\title{Proof of a Conjecture on Overcolored Partition Restricted by Parity of the Parts}
	\author{Imdadul Hussain, Suparno Ghoshal and Arijit Jana}
    \address{Department of Mathematics, National Institute of Technology Silchar, Assam 788010, India}
	\email{imdadulhs05@gmail.com}

	\address{Department of Computer Science, Ruhr University Bochum, Germany}
	\email{suparno.ghoshal@rub.de}
		
    \address{Department of Mathematics, National Institute of Technology Silchar, Assam 788010, India}
    \email{jana94arijit@gmail.com}
	
	
\begin{abstract}
In a recent paper, Thejitha and Fathima introduced the overcolored partition function $\bar{a}_{r,s}(n)$, which enumerates overpartitions in which even parts may appear in one of $r$ colors and odd parts in one of $s$ colors, for fixed integers $r,s \geq 1$. They also proposed several conjectures concerning families of congruences modulo powers of $2$ for specific arithmetic progressions of $\bar{a}_{r,s}(n)$.
In this paper, we provide an elementary proof of this conjecture that relies only on classical $q$-series manipulations and properties of Ramanujan's theta function.
\end{abstract}

\maketitle

\section{Introduction}\label{section1}

We denote by $p_k(n)$ the number of partitions of $n$ into $k$ colors. 
For $q \in \mathbb{C}$ with $|q|<1$, we use the notation
\[
f_n^k := \prod_{j \ge 1} (1 - q^{jn})^{k},
\]
for integers $n,k$ with $n>0$. The generating function for $p_k(n)$ is given by
\begin{align*}
\sum_{n=0}^{\infty} p_k(n) q^n = \frac{1}{f_1^k}.
\end{align*}
In particular, when $k=1$, the function $p_1(n)$ reduces to the unrestricted partition function $p(n)$.

Agarwal and Andrews \cite{agar} were the first to introduce colored partitions in the context of unrestricted partitions. Later, several mathematicians studied different aspects of colored partitions. In particular, Thejitha, Sellers, and Fathima \cite{tsf} considered colored partitions in which the even parts may appear in one of $r$ colors and the odd parts may appear in one of $s$ colors, where $r,s \ge 1$ are fixed integers. They denoted the corresponding counting function by $a_{r,s}(n)$. The generating function for $a_{r,s}(n)$ is given by
\begin{align}\label{gf0}
\sum_{n=0}^{\infty} a_{r,s}(n) q^n=\frac{f_2^{\,s-r}}{f_1^{\,s}}.
\end{align}
The special case $r=1$ was first examined by Hirschhorn and Sellers \cite{hs1}, while the case $s=1$ was  studied by Amdeberhan, Sellers, and Singh \cite{ajsi}.

An overpartition of $n$ is a partition of $n$ in which the first occurrence of each distinct part may be overlined. Let $\overline{p}(n)$ denote the number of overpartitions of $n$, where $\overline{p}(0)=1$. 

The generating function for $\overline{p}(n)$ was introduced by Corteel and Lovejoy \cite{Corteel1} and is given by
\begin{equation*}\label{eqn2}
\sum_{n=0}^{\infty} \overline{p}(n) q^{n}
= \prod_{n=1}^{\infty} \frac{1+q^{n}}{1-q^{n}}
= \frac{f_{2}}{f_{1}^{2}}.
\end{equation*}

A similar generalization appears in the work of Amdeberhan, Sellers, and Singh \cite{ajsi}, where they considered overpartitions in a colored setting. In particular, they extended their earlier study of the case $s=1$ of \eqref{gf0} to overpartitions. They investigated partitions of $n$ in which each even part may occur in one of $r \ge 1$ colors, while the first occurrence of each part may be overlined. Let $\bar{a}_r^*(n)$ denote the number of such partitions.
The generating function for $\bar{a}_r^*(n)$ is given by
\begin{align}\label{1}
\sum_{n=0}^{\infty} \bar{a}_r^*(n) q^n
= \frac{f_4^{\,r-1}}{f_1^2 f_2^{\,2r-3}}.
\end{align}

\noindent Recently, Das, Maity, and Saikia \cite{das} investigated the divisibility properties of $\bar{a}_r^*(n)$ for various values of $r$. In particular, they established the  generalizations of Sellers' results (see \cite[Thms.~2.5 and~2.7]{ovrcu}).

\noindent Very recently, Thejitha and Fathima \cite{ThejithaFathima2026} introduced $\bar{a}_{r,s}(n)$, which denotes the number of partitions of $n$ in which even parts may appear in one of $r$ colors and odd parts may appear in one of $s$ colors, for fixed integers $r,s \ge 1$, with the additional condition that the first occurrence of each part may be overlined. They obtained the generating function for $\bar{a}_{r,s}(n)$ as
\begin{align}\label{gf}
\bar{A}_{r,s}(q):=\sum_{n=0}^{\infty}\bar{a}_{r,s}(n)q^n=\frac{f_2^{\,3s-2r}}{f_1^{\,2s}f_4^{\,s-r}}.
\end{align}
In the same paper, they established several infinite families of congruences for $\bar{a}_{r,s}(n)$ modulo primes and powers of $2$. In particular, their results yield various Ramanujan-type congruences modulo powers of $2$ and modulo primes $p \ge 3$ satisfied by $\bar{a}_{r,s}(n)$. In the concluding remarks, they stated that they would be pleased to see an elementary proof of the congruences proposed in the following conjecture.

\begin{conjecture}\label{conj}\cite[Conjecture 6.1]{ThejithaFathima2026}
    For $k\geq 1$ and $j, i \geq 0$, we have
		\begin{align*}
			\bar{a}_{2^{k+1}j+2^k-1,2i+1}(3n+2)&\equiv 0 \pmod {2^{k+1}},  \\
			\bar{a}_{2^{k+1}j+2^k-1,2i+1}(9n+3)&\equiv 0 \pmod {2^{k+2}}, \\
		      \bar{a}_{2^{k+1}j+2^k-1,2i+1}(9n+6)&\equiv 0 \pmod {2^{k+2}}.
		\end{align*}

\end{conjecture}
The purpose of this paper is to establish an elementary proof of Conjecture~\ref{conj}.
\begin{theorem}
	Conjecture \ref{conj} is true.
\end{theorem}
\section{Proof of Conjecture \ref{conj}}\label{section2}
We have,
\begin{align} \label{p2eq1}
\bar{A}_{r,s}(q) &= \sum_{n=0}^{\infty}\bar{a}_{r,s}(n)q^n = \dfrac{f_2^{3s-2r}}{f_1^{2s}f_4^{s-r}} \notag \\
\text{Put } r = 2^{k+1}j + 2^k - 1, s &= 2^ki + 1 \text{ to obtain} \notag \\
\sum_{n=0}^{\infty}\bar{a}_{2^{k+1}j + 2^k - 1,2^ki + 1}(n)q^n &= \dfrac{f_2^{3\cdot2^ki + 3 - 2^{k+2}j-2^{k+1}+2}}{f_1^{2^{k+1}i+2}f_4^{2^ki + 1 - 2^{k+1}j - 2^k + 1}} \notag \\
&= \dfrac{f_2^{3\cdot2^ki + 3}}{f_1^{2^{k+1}i+2}} \cdot \left(\dfrac{f_4^{2^{k+1}}}{f_2^{2^{k+2}}}\right)^{j} \cdot \dfrac{f_4^{2^{k}}}{f_2^{2^{k+1}}} \cdot \dfrac{f_2^2}{f_4^{2^ki+1}\cdot f_4}
\end{align}
Again, we use the congruence
\[
f_m^{p^k} \equiv f_{mp}^{\,p^{k-1}} \pmod{p^k}.
\]
Applying this with $m=2$ and $p=2$, we obtain
\[
f_2^{2^{k+2}} \equiv f_4^{2^{k+1}} \pmod{2^{k+2}}.
\]

Therefore, from \eqref{p2eq1}, we obtain
\begin{align}\label{p2eq2}
\sum_{n=0}^{\infty}\bar{a}_{2^{k+1}j + 2^k - 1,2^ki + 1}(n)q^n &\equiv  \dfrac{f_2^{3\cdot2^ki + 3}}{f_1^{2^{k+1}i+2}} \cdot \dfrac{f_4^{2^{k}}}{f_2^{2^{k+1}}} \cdot \dfrac{f_2^2}{f_4^{2^ki+1}\cdot f_4} \pmod{2^{k+2}} \notag \\
&= \left(\dfrac{f_2^3}{f_1^2 f_4}\right)^{2^ki+1} \cdot \left(\dfrac{f_4}{f_2^2}\right)^{2^k} \cdot \dfrac{f_2^2}{f_4}\notag\\
&= \left(\dfrac{f_2^5}{f_1^2 f_4^2}\right)^{2^ki+1} \cdot \left(\dfrac{f_4}{f_2^2}\right)^{2^ki} \cdot \left(\dfrac{f_4}{f_2^2}\right)^{2^k}\notag\\
&= \varphi(q)^{2^ki+1} \left(\dfrac{1}{\varphi(-q^2)}\right)^{2^ki} \left(\dfrac{1}{\varphi(-q^2)}\right)^{2^k}\notag\\
&= \varphi(q) \varphi(q)^{2^{k+1}i} \left(\dfrac{1}{\varphi(-q^2)}\right)^{2^{k+1}i} \left(\dfrac{\varphi(-q^2)}{\varphi(q)}\right)^{2^ki} \left(\dfrac{1}{\varphi(-q^2)}\right)^{2^k}\notag\\
&= \varphi(q) \varphi(q)^{2^{k}i} \dfrac{1}{\varphi(-q)^{2^ki}} \left(\dfrac{\varphi(-q)}{\varphi(q)}\right)^{2^{k-1}i} \left(\dfrac{1}{\varphi(q)\varphi(-q)}\right)^{2^{k-1}}\notag\\
&\equiv \varphi(q) \varphi(q)^{2^{k}i} (\varphi(q)\varphi(q^2)^{2})^{2^ki} \left(\dfrac{\varphi(-q)}{\varphi(q)}\right)^{2^{k-1}i} (\varphi(q^2)^2)^{2^{k-1}} \pmod{2^{k+2}} \notag\\
&= \varphi(q) \varphi(q)^{2^{k+1}i} \varphi(q^2)^{2^{k+1}i} \left(\dfrac{\varphi(-q)}{\varphi(q)}\right)^{2^{k-1}i} \varphi(q^2)^{2^{k}}\notag\\
&\equiv \varphi(q) \left(\dfrac{\varphi(-q)}{\varphi(q)}\right)^{2^{k-1}i} \varphi(q^2)^{2^{k}} \pmod{2^{k+2}} 
\end{align}
In the above calculations, we mainly use the identities
\[
\frac{1}{\varphi(-q)}=\varphi(q)\varphi(q^2)^2\varphi(q^4)^4\cdots
\]
and
\[
\varphi(q)\varphi(-q)=\varphi(-q^2)^2.
\]
Moreover, in the final congruence we use the fact that
\[
\varphi(q)^{2^{k+1}i}\equiv 1 \pmod{2^{k+2}}.
\]

Since,
\begin{align}\label{p2eq3}
\dfrac{\varphi(-q)}{\varphi(q)} &= \dfrac{1-2q+2q^4-2q^9+\cdots}{1+2q+2q^4+2q^9+\cdots}\notag\\
&= (1-2q+2q^4-2q^9+\cdots)(1-2q+4q^2-8q^3+\cdots)\notag\\
&= (1 - 4q + 8q^2 - 16q^3 + 32q^4 - \cdots)\notag\\
&= 1 + 4(-q + 2q^2 - 4q^3 + 8q^4 + \cdots)\notag\\
&= 1 + 4H(q)\\
\text{where } H(q) &= -q + 2q^2 - 4q^3 + 8q^4 + \cdots.  \notag
\end{align}

And, 
\begin{align}\label{p2eq4}
\varphi(q^2)^{2^{k}} &= (1 + 2\sum\limits_{n \geq 1} q^{2n^2})^{2^k}\notag\\
&= 1 + 2^{k+1}\sum\limits_{n \geq 1} q^{2n^2} + \binom{2^k}{2} \left(2\sum\limits_{n \geq 1} q^{2n^2}\right)^2 + \sum\limits_{m=3}^{2^k} \binom{2^k}{m}\left(2\sum\limits_{n \geq 1} q^{2n^2}\right)^m\notag \\
&= 1 + 2^{k+1}\sum\limits_{n \geq 1} q^{2n^2} + 2^{k+1}(2^k-1) \left(\sum\limits_{\alpha, \beta \geq 1} q^{2(\alpha^2 + \beta^2)}\right) + \sum\limits_{m=3}^{2^k} \binom{2^k}{m}\left(2\sum\limits_{n \geq 1} q^{2n^2}\right)^m 
\end{align}

We can write 
\begin{align*}
\binom{2^k}{m} &= \dfrac{2^k(2^k-1)(2^k-2) \cdots (2^k-m+1)}{m!}= \dfrac{\prod\limits_{j=0}^{m-1}(2^k-j)}{m!}
\end{align*}
Let $\nu_p(a)$ denote the $p$-adic valuation of $a$, that is, the highest power of $p$ that divides $a$.
We know, 
\begin{align*}
\nu_2\left(\prod\limits_{j=0}^{m-1}(2^k-j)\right) 
&= k + \sum\limits_{j=1}^{m-1}\nu_2\left(2^k-j\right) = k + \sum\limits_{j=1}^{m-1}\nu_2\left(j\right)
\end{align*} 
and 
\begin{align*}
\nu_2(m!) &= \sum\limits_{j=1}^m \nu_2(j).
\end{align*}
Thus, 
\begin{align*}
\nu_2\binom{2^k}{m} &= k + \sum\limits_{j=1}^{m-1}\nu_2\left(j\right) - \sum\limits_{j=1}^m \nu_2(j)= k - \nu_2(m)
\end{align*}
Hence, $\nu_2\left(2^m\binom{2^k}{m}\right) = m + k - \nu_2(m)$.
For $m \geq 3$, $m + k - \nu_2(m) \geq k+2$.
Therefore, $2^m\binom{2^k}{m} \equiv 0 \pmod{2^{k+2}}$ for $m\geq 3$.
This gives, 
\begin{align}\label{p2eq5}
\varphi(q^2)^{2^{k}} &\equiv 1 + 2^{k+1}\sum\limits_{n \geq 1} q^{2n^2} + 2^{k+1} \left(\sum\limits_{\alpha, \beta \geq 1} q^{2(\alpha^2 + \beta^2)}\right) \pmod{2^{k+2}}\notag \\
&= 1 + 2^{k+1}F(q)\\
\text{where } F(q) &= \sum\limits_{n \geq 1} q^{2n^2} +  \left(\sum\limits_{\alpha, \beta \geq 1} q^{2(\alpha^2 + \beta^2)}\right). \notag
\end{align}
Combining \eqref{p2eq2}, \eqref{p2eq3}, and \eqref{p2eq5}, we obtain
\begin{align}\label{p2eq6}
\sum_{n=0}^{\infty}\bar{a}_{2^{k+1}j + 2^k - 1,2^ki + 1}(n)q^n &\equiv
\varphi(q) \left(\dfrac{\varphi(-q)}{\varphi(q)}\right)^{2^{k-1}i} \varphi(q^2)^{2^{k}} \notag \\
&= \varphi(q)\left[1 + 4H(q)\right]^{2^{k-1}i} [1 + 2^{k+1}F(q)]\notag \\
&\equiv \varphi(q)[1 + 2^{k+1}iH(q)] [1 + 2^{k+1}F(q)] \pmod{2^{k+2}}\notag\\
&\equiv ( 1 + 2\sum\limits_{n \geq 1} q^{n^2}) [1 - 2^{k+1}i q] [1 + 2^{k+1}F(q)] \pmod{2^{k+2}}\notag\\
&\equiv ( 1 + 2\sum\limits_{n \geq 1} q^{n^2}) [1 - 2^{k+1}i q + 2^{k+1}F(q)] \pmod{2^{k+2}}\notag\\
&\equiv  1 + 2\sum\limits_{n \geq 1} q^{n^2} - 2^{k+1}i q + 2^{k+1}F(q) \pmod{2^{k+2}}\notag\\
&\equiv 1 + 2\sum\limits_{n \geq 1} q^{n^2} - 2^{k+1}i q + 2^{k+1}\left[\sum\limits_{n \geq 1} q^{2n^2} - \sum\limits_{\alpha, \beta \geq 1} q^{2(\alpha^2 + \beta^2)}\right] \pmod{2^{k+2}}.
\end{align}
The third congruence holds because of the fact $H(q) \equiv -q \pmod{2}.$
The possible residues of $n^2 \pmod{9}$ are $\{0,1,4,7\}$. So, the coefficient of $q^{9n+3}$ and $q^{9n+6}$ is 0 in 
$\sum\limits_{n \geq 1} q^{n^2}$. The possible residues of $2n^2 \pmod{9}$ are $\{0,2,5,8\}$. So, the coefficient of $q^{9n+3}$ and $q^{9n+6}$ is 0 in $\sum\limits_{n \geq 1} q^{2n^2}$. The possible residues of $2(\alpha^2+\beta^2) \pmod{9}$ are $\{0,1,2,4,5,7,8\}$. So, the coefficient of $q^{9n+3}$ and $q^{9n+6}$ is 0 in 
$\sum\limits_{\alpha, \beta \geq 1} q^{2(\alpha^2 + \beta^2)}$.

As a result, 
\begin{align*}
\bar{a}_{2^{k+1}j + 2^k - 1,2^ki + 1}(9n+3) &\equiv 0 \pmod{2^{k+2}},\\
\bar{a}_{2^{k+1}j + 2^k - 1,2^ki + 1}(9n+6) &\equiv 0 \pmod{2^{k+2}}. \tag*{\qed}
\end{align*}

Also in modulo $2^{k+1}$ \eqref{p2eq6} gives us 
\begin{align*}\label{p2eq7}
\sum_{n=0}^{\infty}\bar{a}_{2^{k+1}j + 2^k - 1,2^ki + 1}(n)q^n &\equiv  1 + 2\sum\limits_{n \geq 1} q^{n^2} \pmod{2^{k+1}}.
\end{align*}
The possible residues of $n^2 \pmod{3}$ are $\{0,1\}$. So, the coefficient of $q^{3n+2}$ is 0 in 
$\sum\limits_{n \geq 1} q^{n^2}$.
Therefore,
\begin{align*}
\bar{a}_{2^{k+1}j + 2^k - 1,\,2^k i + 1}(3n+2) \equiv 0 \pmod{2^{k+1}}. \tag*{\qed}
\end{align*}
This completes the proof.

    {\bf Author contributions:} All authors have contributed equally to the preparation of this manuscript.
	
	{\bf Conflicts of Interest:} The authors do not have any conflict
	of interest.
    
    {\bf Data Availability Statement:}
	No datasets were generated or analyzed in this study. The results are obtained entirely through theoretical and mathematical methods.

\end{document}